# Some addenda on distance function wavelets

W. Chen

Simula Research Laboratory, P. O. Box. 134, 1325 Lysaker, Norway

E-mail: wenc@simula.no

(9 July 2002)

**Summary**

This report will add some supplements to the recently finished report series on the distance function wavelets (DFW). First, we define the general distance in terms of the Riesz potential, and then, the distance function Abel wavelets are derived via the fractional integral and Laplacian. Second, the DFW Weyl transform is found to be a shifted Laplace potential DFW. The DFW Radon transform is also presented. Third, we present a conjecture on truncation error formula of the multiple reciprocity Laplace DFW series and discuss its error distributions in terms of node density distributions. Forth, we point out that the Hermite distance function interpolation can be used to replace overlapping in the domain decomposition in order to produce sparse matrix. Fifth, the shape parameter is explained as a virtual extra axis contribution in terms of the MQ-type Possion kernel. The report is concluded with some remarks on a range of other issues.

***Keywords***: *distance function wavelets, partial differential equation, eigenfunction, general solution, fundamental solution, differentiation elimination, general distance, Abel transform, fractional integral, fractional Laplacian,, Weyl transform, Radon transform, domain decomposition, edge effect, shape parameter, truncation error, error distribution, compactly supported DFW, geodesic MQ-type distance function.*

# 1. Introduction

This report presents supplementary results to recent report series [1-3] on the distance function wavelets (DFW) transforms and series. The classic integral transforms and series are seen to base on eigenfunctions of partial differential equation (PDE), while the DFW uses the fundamental and/or general solutions of PDEs. Notably different from the former two, the common wavelets have very little to do with the solutions of PDEs and are often constructed via the filter theory and spline construction approach. Similarly, most research on the distance functions (radial basis functions, rbf) has long ignored the underlying PDE and integral equation theory. The PDE is arguable language describing the laws of universe [4]. Thus, the wavelets and distance functions should be established on the kernel solutions of PDEs. [5] has made a brief comment on wavelets (rbf prewavelets) and rbf from a PDE perspective.

The standard applications of the distance functions such as the multiquadratics (MQ) and thin plate spline (TPS) are simply single-scale calculations lacking crucial multiscale analysis capability. By now, many efforts have been dedicated to examining the invertiblility of the distance function approximations. Notably, scanty research, however, has involved with issues of completeness and convergence. In the author's view, the DFW series provides a complete basis for distance function approximations and thus guarantees the convergence. Most importantly, the DFW are very powerful and flexible to deal with the multifractal, multiscale, multivariate scattered data and meshfree numerical solution of PDE.

The rest of the report is organized into a few independent sections. Section 2 further refines the general (fractional) distance in terms of the Riesz potential. The Abel DFW transform is also presented via the fractional integral and Laplacian. Section 3 illustrates the underlying connections between the DFW Laplace transform, Weyl transform and Radon transform. In section 4, by analogy with the Taylor and polynomial interpolation expansions, a conjecture is presented on the truncation error formula of the multiple reciprocity (MR) Laplace DFW series, and then, the error distribution is discussed in comparison to the Chebyshev polynomial interpolation. In section 5, the local Hermite distance function expansion is presented to compete with overlapping in the domain decomposition. In terms of the Possion kernel, section 6 defines the MQ shape parameter as an extra axis contribution to Euclidean distance. Finally, some issues are remarked in section 7.

## 2. General distance, fractal, fractional integral and Laplacian, Abel DFW transforms

The Riesz potential [6] is defined by

$$I_s f(x) = \frac{\Gamma((n-s)/2)}{\pi^{n/2} 2^s \Gamma(s/2)} \int_{IR^n} \frac{f(\xi)}{\|x-\xi\|^{n-s}} d\xi, \qquad 0 \prec s \prec n, \qquad (1)$$

where $s$ is the order of Riesz potential, and $n$ is the topological dimension. Accordingly, [3] defined a general distance variable by

$$r_k^s = \left( \sum_{i=1}^n (x^i - x_k^i)^s \right)^{1/s}, \qquad (2)$$

where $s$ could even be a complex number in relation to complex dimensions and is connected with capacity dimension, fractional derivative, fractal [3] and Holder exponent. It is noted that $s=2$ leads to the Newtonian potential and corresponding Euclidean distance.

It is known [7] that the forward Abel transform

$$f(x) = \int_0^x \frac{g(t)dt}{(x-t)^{1-\beta}} \qquad 0 \prec \beta \prec 1, \qquad (3)$$

can be seen as a fractional integral [8] in one dimension, while the backward Abel transform

$$g(t) = -\frac{\sin(\pi(1-\beta))}{\pi} \frac{d}{dt} \int_0^t \frac{f(x)dx}{(x-t)^\beta} \tag{4}$$

is easily regarded as a fractional derivative. [6] also relates the fractional integral to the Riesz potential [1]. Thus, we can construct the multidimensional DFW Abel transforms

$$H(y,x) = \frac{\Gamma((n-y)/2)}{\pi^{n/2} 2^y \Gamma(y/2)} \int_{IR^n} \frac{h(\xi)}{\|x-\xi\|^{n-y}} d\xi, \quad 0 \prec y \prec 2, \tag{5}$$

$$h(x) = \left(-\Delta^{y/2}\right) H(y,x), \tag{6}$$

where $-\Delta^{y/2}$ is the fractional Laplacian of $y$ order.

## 3. DFW Laplace transforms and Radon transforms

The DFW Weyl transform was presented by [3]:

$$W_n(\xi,\gamma) = \frac{2\Gamma(n/2)}{\sqrt{\pi}\Gamma\left(\frac{n-1}{2}\right)} \int_{IR^n} f(x) \left[\|\xi-x\|^2 - \gamma^2\right]^{(n-3)/2} \|\xi-x\| H\left(\|\xi-x\|^2 - \gamma^2\right) dx, \tag{7}$$

where $H$ is the Heaviside step function. [9] links the radial Weyl and Radon transforms by

$$R_n f = \frac{\pi^{n/2}}{\Gamma(n/2)} W_n f \,. \qquad (8)$$

By (7) and (8), we have the DFW Radon transform

$$R_n(\xi,\gamma) = \frac{2\pi^{n/2}}{\Gamma(n/2)} \int_{IR^n} f(x) \left[\|\xi - x\|^2 - \gamma^2\right]^{(n-3)/2} \|\xi - x\| H\left(\|\xi - x\|^2 - \gamma^2\right) dx, \qquad (9)$$

Let the shape parameter $\gamma=0$, the DFW Radon transform turns out to be the DFW Laplace potential transform

$$\begin{aligned} D_n(\xi,n) &= \frac{2\pi^{n/2}}{\Gamma(n/2)} \int_{IR^n} f(x) \|\xi - x\|^{(n-2)} dx \\ &= \frac{1}{n-2} \int_{IR^n} \frac{f(x)}{u^*_{L_n}(\|\xi - x\|)} dx \end{aligned} \qquad n\neq 2, \qquad (10)$$

where $u^*_{L_n}$ is the fundamental solution of the $n$-dimensional Laplace equation, and the "scale parameter" of this DFW transform is the dimensionality. It is easy to see that the Weyl transform is a shifted DFW Laplace transform.

According to (10), we find that the weight function $1/\|\xi - x\|$ in various potential DFW transforms presented in chapter 2 and 3 of [3] may not be necessary. In addition, the dimensional parameter $n$ in the backward DFW potential transforms may be replaced by (n-2) instead or simply removed.

## 4. Conjecture on truncation error of MR Laplace DFW series

The multiple reciprocity Laplacian DFW series is given by [3]

$$f(x) = f_0(x) + \sum_{m=1}^{N}\sum_{k=1}^{M} \alpha_{mk} u^*_{L_n^m}(\|x - x_k\|) + R(x), \qquad (11)$$

where $f_0(x)$ is evaluated via the boundary integral equation detailed in [1,3], and $u^*_{L_n^m}$ is the fundamental solution of the $m$-th order n-dimensional Laplacian $\nabla^{2(m+1)}$. $R(x)$ is the truncation error. It is noted that when the order $m \geq 1$, these kernel solutions are no longer singular at the origin.

The above MR Laplacian DFW series is an analogue to the polynomial interpolation. By analogy with the truncation error formula of Taylor and polynomial interpolation expansions, a conjecture formula on the MR Laplacian DFW expansion (11) is given by

$$R(x) = \frac{C}{(M!)^N} \prod_{k}^{M} \|x - x_k\| \left(\nabla^{2(N+1)} f(\xi)\right), \qquad (12)$$

where $\xi$ is a node among interpolation domain. [10] conjectured that the truncation error under $M$ nodes and $N$ scale might be $O(M^N (\log M)^{n-1})$ and is independent of dimensionality. The known accuracy orders of the TPS and cubic spline are the special cases of (12).

As in the case of polynomial interpolations, the error distribution is determined by

$$w(x) = \prod_{k}^{M} \|x - x_k\|. \tag{13}$$

The optimal placement of nodes should satisfy the minimization of maximum $w(x)$, namely,

$$x_k \to Min(Max(w(x))). \tag{14}$$

In the classic Chebyshev polynomials, we have $L_2$ optimal distribution of truncation errors which satisfies (14). Under the uniformly-spaced nodes, the truncation errors at the central region will be much less than those at the vicinity of the boundary. So, the nodes should be placed increasingly dense toward the boundary as the zeros of Chebyshev polynomials in one dimensional case. Meanwhile, the condition of interpolation matrix under the optimal placement of nodes is also relatively better.

### 5. Hermite interpolation, edge effect, and domain decomposition

The domain decomposition seems a promising technique to result in a sparse matrix of distance function interpolation. Overlapping among different subdomains is a common approach for connectivity. However, the arbitrary placement of overlapping nodes causes some troublesome issues in practical engineering computing.

Let $Q(x)$ is the object function to be approximated. If

$$\Re\{Q(x)\} = g(x) \tag{15}$$

where $\mathscr{R}$ is a differential operator, and $g(x)$ tends to be zero or very smooth. We think $\mathscr{R}$ is a good smoothing operator in the distance function approximation of $Q(x)$. In terms of Green second identity, we have

$$c(x)Q(x) = -\int_\Omega g(x_k)u_\mathscr{R}^*(x-x_k)dx_k + \int_{\partial\Omega}\left[\left\{\frac{\partial Q(x_j)}{\partial n}u_\mathscr{R}^*(x-x_j) - Q(x_j)\frac{\partial u_\mathscr{R}^*(x-x_j)}{\partial n}\right\}\right]dx_j \quad (16)$$

where $u_\mathscr{R}^*$ is the fundamental solution of $\mathscr{R}$; $n$ is the unit outward normal; $x_k$ and $x_j$ denote source points on the whole domain and the boundary, respectively. $c(x)$ depends on the location, $c(x)=1$ inside domain and $c(x)=1/2$ on smooth boundary, varying $c(x)$ on irregular boundary vertices.

By numerical integration of (16), Chen [11] presented a Hermite scheme of distance function interpolation to mitigate the so-called edge effect in the function interpolation and numerical solution of PDEs, i.e.

$$Q(x) = \sum_{i=1}^{N+L}\beta_i\varphi(x-x_i) + \sum_{i=N+L+1}^{N+2L}\beta_i\varphi(x-x_i) + \sum_{i=N+L+1}^{N+3L}\beta_i\left(-\frac{\partial\varphi(x-x_i)}{\partial n}\right), \quad (17)$$

where $\varphi$ is the kernel distance function built on the fundamental solution $u_\mathscr{R}^*$. $N$ and $L$ are respectively the numbers of inside domain and boundary nodes. Note that the boundary nodes are here interpolated twice respectively for Dirichlet and Neumann data. To keep the symmetry of the interpolation matrix, a minus sign is added before Neumann boundary in (17) due to the fact that Neumann boundary condition is not self-adjoint.

In terms of the Local boundary integral equation [12] and corresponding Green second identity, we can extend the above interpolation scheme into subdomains and has an alternative scheme to overlapping in the domain decomposition technique. Namely, by keeping the different degrees of differential continuity, we can connect different subdomains without using overlapping nodes. The penalty is that the scheme is not fully meshfree any more since we need connect large or small subregions through somehow connectivity like those in the standard FEM mesh.

(17) can also be easily extended to implementing multiple boundary conditions of high order partial equation problems. For instance, we need to enforce two boundary conditions at each boundary knot for the biharmonic equation. In that case, we have

$$Q(x) = \sum_{i=1}^{N+L} \beta_i \varphi(x - x_i) + \sum_{i=N+1}^{N+L_1} \beta_i B_1\{\varphi(x - x_i)\} + \sum_{i=N+L_1+1}^{N+L_1+L_2} \beta_i B_2\{\varphi(x - x_i)\} + \sum_{i=N+L_1+L_2+1}^{N+3L} \beta_i B_3\{\varphi(x - x_i)\},$$

(18)

where B's and L's with the subscripts 1,2,3 indicate the three different boundary conditions and the corresponding numbers of boundary knots. For Neumann boundary conditions we need add a minus sign as did in (15) to get a symmetric interpolation.

**6. Shape parameter in MQ-type kernel distance functions**

Replacing the Euclidean distance variable $r_k$ by the MQ $\sqrt{r_k^2 + c^2}$ in the rotational invariant fundamental solution and general solution of PDEs will yield a variety of the MQ-type kernel distance functions (KDF) [3,10]. In particular, [3] pointed out that the Possion kernel belongs to the MQ-type KDF, i.e.

$$P_q^n(\|x - x_k\|, q) = \frac{c_n q}{\left(\|x - x_k\|^2 + q^2\right)^{(n+1)/2}}, \qquad q>0, \tag{19}$$

where $n$ denotes dimension, and $x \in R^{n+1}$. $\int P_q^n(x)dx = 1$ for all $q>0$, and $q$ is the shape parameter in terms of the MQ. The object function $f(x)$ is defined boundary condition on $q=0$ surface of a $R^{n+1}$ domain, i.e.

$$u(x,0) = f(x). \tag{20}$$

We also have

$$u(x,q) = (P_q^n * f)(x) = \frac{2q}{\omega_{n+1}} \int_{IR^n} \frac{f(\xi)}{\left(\|\xi - x\|^2 + q^2\right)^{(n+1)/2}} d\xi. \tag{21}$$

When shape parameter $q \to \infty$, the function is approximated as the boundary value in infinite surface.

The MQ-type distance function expansion should be stated as

$$f(x) \cong a_0 + \sum_{j=1}^{N} \sum_{k=1}^{M} \alpha_{mk} \left(\sqrt{\|x - x_k\| + c_j^2}\right). \tag{22}$$

It is noted that the augmented polynomial terms are not necessary as in the boundary integral equation. We also have a general DFW series as follows:

$$f(x) \cong a_0 + \sum_{m=1}^{N} \sum_{k=1}^{M} \alpha_{mk} \|x - x_k\|^m. \tag{23}$$

## 7. Other issues

The compactly-supported rbfs (CS-RBF) [13,14] are built with the spline construction. [15] also creates the compactly-supported spline rbf wavelets. Both CS-RBF [13,14] and CS-RBF wavelets [15] do not consider the dimensional effect. By analogy with the compactly-supported polynomial basis function in the FEM, we can also construct the CS DFW via the kernel solution of the Laplace equation, where the building blocks will be the fundamental solutions of different order Laplacians. We can call them the Laplacian spline, which also has something to do with what we discussed in the preceding section [5].

Furthermore, by analogy with the rational polynomial interpolation, we also have the rational Laplace DFW series of the form:

$$P(x) = \frac{P_0(x) + \sum_{m=1}^{N_1}\sum_{k=1}^{M_1} \alpha_{mk} u^*_{L_n^m}(\|x - x_k\|)}{g_0(x) + \sum_{i=1}^{N_2}\sum_{j=1}^{M_2} \beta_{ij} u^*_{L_n^i}(\|x - x_j\|)} \quad . \tag{24}$$

The DFW built on the fundamental and general solutions of PDEs are scale orthogonal and can be orthonormalized via the Gram-Schmidt orthonormalization scheme over irregular translation. For regular translation, it is not difficult to get the scale and translation orthogonal DFWs.

The Gaussian rbf is often considered as a candidate of weight function in moving least squares in relation to various meshfree methods. In fact, the kernel distance function presented in [1-3] can be used to serve the same function. The differentiation smoothing defined in [3] may be better renamed as differentiation elimination, which eliminates the inhomogeneous terms via the composite differentiation processing and leads to a boundary-only approximation scheme. We also need to mention that the Neumann and Dirichlet eigenvalues of HF series discussed in section 3.3 of [1] correspond to the respective homogeneous boundary conditions and have nothing to do with the known inhomogeneous boundary conditions.

By analogy with the high-order fundamental solutions of the geodesic Laplacian distance functions [3], we have the corresponding geodesic TPS and cubic splines for heterogeneous media:

$$\phi(x-\xi) = \begin{cases} [\kappa_{ij}]^{-1/2} R^{2m} \ln R, & n=2, \\ [\kappa_{ij}]^{-1/2} R^{2m-1}, & n=3, \end{cases} \quad (25)$$

where the coefficient matrix $\kappa = [\kappa_{ij}]$ represents the parameters in different directions ($j$) and locations ($i$) of anisotropic and inhomogeneous media, and the geodesic distance $R$ is defined by

$$R^2 = \sum_{i,j=1}^{n} \kappa_{ij}^{-1}(x_i - \xi_i)(x_j - \xi_j). \quad (26)$$

Similarly, we have the geodesic MQ-type (Gaussian) distance functions:

$$\varphi(x-\xi) = [\kappa_{ij}]^{-1/2} \sqrt{R^2 + c^2}, \quad (27a)$$

$$g(x-\xi) = \frac{[\kappa_{ij}]^{-1/2}}{\sqrt{R^2 + c^2}}, \quad (27b)$$

$$\psi(x-\xi) = \begin{cases} [\kappa_{ij}]^{-1/2} (R^2 + c^2)^m \ln(R^2 + c^2), & n=2, \\ [\kappa_{ij}]^{-1/2} (R^2 + c^2)^{2m-1}, & n=3, \end{cases} \quad (28)$$

$$P_c^n(x-\xi) = \frac{[\kappa_{ij}]^{-1/2}}{(R^2+c^2)^{(n+1)/2}}, \quad (29)$$

$$G(x-\xi) = [\kappa_{ij}]^{-1/2} e^{-R^2/\alpha^2}, \quad (30)$$

where $c$ and $\alpha$ are the shape parameter, and $n$ is the dimensionality.